%% file: athn_cp.tex
\documentclass[a4paper,USenglish,cleveref, autoref, thm-restate]{lipics-v2021}

\usepackage[short,nocomma]{optidef}
\usepackage{tikz}
\usetikzlibrary{calc,arrows,positioning}
\usepackage{booktabs}
\usepackage{pifont}

\newcommand{\duration}{\textrm{duration}}
\newcommand{\type}{\textrm{type}}
\newcommand{\hub}{\textrm{hub}}
\newcommand{\task}{\textrm{task}}
\newcommand{\load}{\texttt{l}}
\newcommand{\unload}{\texttt{u}}
\newcommand{\drive}{\texttt{d}}
\newcommand{\relocate}{\texttt{r}}
\newcommand{\park}{\texttt{p}}
\definecolor{loadcolor}{HTML}{f4a582}
\definecolor{unloadcolor}{HTML}{f4a582}
\definecolor{relocatecolor}{HTML}{ca0020}
\definecolor{parkcolor}{HTML}{92c5de}
\definecolor{drivecolor}{HTML}{0571b0}
\newcommand{\cmark}{\ding{51}}
\newcommand{\xmark}{\ding{55}}

\pdfoutput=1 
\hideLIPIcs  

\graphicspath{{./figures/}}

\title{Constraint Programming to Improve Hub Utilization in Autonomous Transfer Hub Networks} 

\titlerunning{Constraint Programming to Improve Hub Utilization in ATHNs} 

\author{Chungjae Lee}{H. Milton Stewart School of Ind. and Syst. Engineering, Georgia Institute of Technology, USA \and \url{https://www.isye.gatech.edu/users/chungjae-lee}}{clee384@gatech.edu}{https://orcid.org/0000-0002-9857-1789}{}
\author{Wirattawut Boonbandansook}{H. Milton Stewart School of Ind. and Syst. Engineering, Georgia Institute of Technology, USA}{maxboonban@gatech.edu}{https://orcid.org/0009-0004-0306-1113}{}
\author{Vahid Eghbal Akhlaghi}{H. Milton Stewart School of Ind. and Syst. Engineering, Georgia Institute of Technology, USA \and \url{https://www.isye.gatech.edu/users/vahid-eghbal-akhlaghi}}{vahid.eghbal@gatech.edu}{https://orcid.org/0000-0002-3120-4108}{}
\author{Kevin Dalmeijer\footnote{Corresponding author}}{H. Milton Stewart School of Ind. and Syst. Engineering, Georgia Institute of Technology, USA \and \url{https://www.isye.gatech.edu/users/kevin-dalmeijer}}{dalmeijer@gatech.edu}{https://orcid.org/0000-0002-4304-7517}{}
\author{Pascal Van Hentenryck}{H. Milton Stewart School of Ind. and Syst. Engineering, Georgia Institute of Technology, USA \and \url{https://www.isye.gatech.edu/users/pascal-van-hentenryck}}{pvh@gatech.edu}{https://orcid.org/0000-0001-7085-9994}{}

\authorrunning{C. Lee, W. Boonbandansook, V. Eghbal Akhlaghi, K. Dalmeijer, P. Van Hentenryck} 

\Copyright{Chungjae Lee, Wirattawut Boonbandansook, Vahid Eghbal Akhlaghi, Kevin Dalmeijer, Pascal Van Hentenryck} 



\category{Short Paper} 

\relatedversion{} 

\acknowledgements{
	This research was partly funded through a gift from Ryder and partly supported by the NSF AI Institute for Advances in Optimization (Award 2112533). Special thanks to the Ryder team for their invaluable support, expertise, and insights.
}

\nolinenumbers 

\begin{document}

\maketitle

\begin{abstract}
The Autonomous Transfer Hub Network (ATHN) is one of the most promising ways to adapt self-driving trucks for the freight industry.
These networks use autonomous trucks for the middle mile, while human drivers perform the first and last miles.
This paper extends previous work on optimizing ATHN operations by including transfer hub capacities, which are crucial for labor planning and policy design.
It presents a Constraint Programming (CP) model that shifts an initial schedule produced by a Mixed Integer Program to minimize the hub capacities.
The scalability of the CP model is demonstrated on a case study at the scale of the United States, based on data provided by Ryder System, Inc.
The CP model efficiently finds optimal solutions and lowers the necessary total hub capacity by 42\%, saving \$15.2M in annual labor costs.
The results also show that the reduced capacity is close to a theoretical (optimistic) lower bound.
\end{abstract}

\input{sections/introduction.tex}
\input{sections/methodology.tex}
\input{sections/case_study.tex}
\input{sections/results.tex}
\input{sections/conclusion.tex}



\bibliography{athn_cp}

\end{document}

%% file: sections/introduction.tex
\section{Introduction}
\label{sec:introduction}

It is widely believed that autonomous trucks will revolutionize the freight industry, and many companies have started exploring its potential \cite{FleetOwner2021-TusimpleAutonomousTruck, Forbes2021-PlusPartnersIveco, FreightWaves2021-GatikIsuzuPartner, HDT2021-DaimlersRedundantChassis, TechCrunch2021-AuroraVolvoPartner, TransportTopics2020-NavistarTusimplePartner}.
Some of the major players describe the \emph{transfer hub business model} to be the most likely implementation for autonomous trucking \cite{RolandBerger2018-ShiftingGearAutomation, ShahandashtEtAl2019-AutonomousVehiclesFreight, Viscelli-Driverless?AutonomousTrucks}.
An Autonomous Transfer Hub Network (ATHN) is a network of autonomous truck ports (\emph{transfer hubs}) that leverages the strengths of humans and automation in their most effective roles.
Autonomous trucks handle the monotonous middle mile to transport goods between the transfer hubs, while humans handle the complex first and last miles through local cities and deal with customer contacts.
According to Roland Berger \cite{RolandBerger2018-ShiftingGearAutomation}, implementing a transfer hub model can lower operational costs by 22\% to 40\%.
A case study in the Southeast of the United States by Ryder System, Inc. and the Socially Aware Mobility lab \cite{RyderSAM2021-ImpactAutonomousTrucking} reports savings from 27\% to 40\%, supporting earlier estimates.
The authors model optimizing ATHN operations as a scheduling problem, and a Constraint Programming (CP) model is used to minimize the empty miles \cite{DalmeijerVanHentenryck2021-OptimizingFreightOperations}.
Subsequent work by \cite{LeeEtAl2022-OptimizationModelsAutonomous} presents a column-generation approach and a bespoke network-flow model to quickly find high-quality solutions.
These findings are combined by \cite{LeeEtAl2023-OptimizingAutonomousTransfer} into a flow-based Mixed Integer Programming (MIP) framework that can solve large-scale instances over a long time horizon (e.g., a month) optimally in reasonable time. These capabilities also allow for detailed analyses of the benefits and costs of ATHNs.

This paper builds on previous work by introducing a framework for optimizing both ATHN operations and transfer hub capacity utilization.
Capacity planning plays a critical role in an actual implementation of the ATHN for operational scheduling.
Furthermore, hub capacity has an impact on the number of essential personnel required, which is a crucial factor to consider in labor planning and policy design.
Prior research on ATHN operations primarily focused on routing and scheduling, but neglected the capacity considerations involved.
Including capacity constraints in a scalable way is not obvious.
The CP method in \cite{DalmeijerVanHentenryck2021-OptimizingFreightOperations} could be extended to incorporate hub capacity constraints by modeling a Resource-Constrained Project Scheduling Problem \cite{HARTMANN20221, HERROELEN1998279, Pritsker1969MultiprojectSW} at every hub and adding cumulative constraints \cite{AGGOUN199357c, Schutt2011} for hub capacity.
However, the poor performance reported in \cite{DalmeijerVanHentenryck2021-OptimizingFreightOperations} makes it unlikely this method will provide good solutions on the national level.
The MIP framework from \cite{LeeEtAl2023-OptimizingAutonomousTransfer} does handle large scale systems, but does not support a cumulative constraint.
Alternatives typically require sophisticated modeling techniques and solution methods such as Time-expanded Networks or Dynamic Discretization Discovery \cite{Boland_CSND, Perspectives_tdmodels, He_exact_sndhc, LukeMarshall_Interval_DDD, ScherrEtAl2020-DynamicDiscretizationDiscovery, DDD_TDTSPTW, Wu_SNDHC}, which are not obvious to scale either.

The method presented in this work optimizes ATHN operations and improves hub utilization even for large-scale systems by combining the strengths of MIP and CP.
The MIP model is used to generate routes and an initial schedule, while the CP model is used to shift the schedule and minimize the required hub capacities.
A case study based on real data is conducted to demonstrate the effectiveness of the new methodology on an ATHN system spanning the United States for a four-week horizon.
The proposed CP model efficiently finds optimal solutions and lowers the necessary total hub capacity by 42\%.
This reduction in capacity may save \$15.2M per year in labor cost.
Furthermore, it is shown that this is close to the best possible savings for any initial schedule.
This paper also includes a sensitivity analysis and provides operational insights for future implementation of the ATHN framework.
The remainder of this paper is organized as follows.
Section~\ref{sec:methodology} presents the MIP and CP methodology.
Section~\ref{sec:case_study} describes the data and experimental settings used in the case study.
Results and sensitivity analysis are presented in Section~\ref{sec:results}.
Finally, Section~\ref{sec:conclusion} provides conclusions and suggests directions for future research.

%% file: sections/methodology.tex
\section{Methodology}
\label{sec:methodology}

This paper uses the methods by \cite{LeeEtAl2023-OptimizingAutonomousTransfer} to design an ATHN and to optimize the routes of the autonomous vehicles with a MIP.
The resulting solution minimizes the cost of the system, but does not take into account the necessary capacity at each of the hubs, which may lead to low hub utilization.
To address this issue, this paper introduces a CP model to shift the schedule in such a way that the original time windows remain satisfied, and the necessary hub capacities are minimized.

\subsection{ATHN Design and Operations}
The input data for the design and optimization of ATHN consists of a set of loads with origins, destinations, and release times.
Following \cite{LeeEtAl2023-OptimizingAutonomousTransfer}, K-means clustering is used to determine hub locations $H$, and loads are assigned to hubs according to a hub-assignment rule that minimizes the total driving distance (with autonomous miles discounted by a factor $\gamma$).

To optimize the ATHN operations, first define a set of tasks $T = \{1, 2, \hdots \}$, where each task $t \in T$ corresponds to moving a load from origin hub $h^+_t \in H$ to destination hub $h^-_t \in H$ with an autonomous truck.
Each task $t \in T$ is associated with a desired pickup time $p_t$ and a flexibility $\Delta \ge 0$, resulting in a pickup time window of $[p_t - \Delta, p_t + \Delta]$.
A \emph{task graph} $G=(V,A)$ is introduced to find the optimal sequence of tasks for every vehicle.
The set $V = T \cup \{0\} \cup \{\lvert T \rvert + 1\}$ consists of vertices that correspond to the tasks, together with a source node $0$ and a sink node $\lvert T \rvert + 1$.
Choosing an arc $a\in A$ indicates that the corresponding tasks are performed sequentially by the same vehicle.
An arc between two tasks $t, t' \in T$ represents loading at $h^+_t$, moving freight from $h^+_t$ to $h^-_t$, unloading at $h^-_t$, and relocating from $h^-_t$ to $h^+_{t'}$ to be ready for the next task.
Each arc is associated with a corresponding time $\tau_a$.
The time for loading or unloading is given by a parameter $\sigma$, and OpenStreetMap times are used for driving and relocation \cite{OpenStreetMap2021}.

The cost is calculated in two parts: $d_t$ is the \emph{direct cost}, which represents the cost of serving task $t\in T$ directly with a conventional truck without using any of the hubs.
This cost is taken to be the total distance for delivery and empty return.
Note that the current non-autonomous system corresponds to using only direct trips.
The second component is a \emph{cost differential} $c_a$ associated with each arc $a \in A$.
For arc $a = (t, t')$ this represents the difference in cost to switch task $t$ from conventional to autonomous delivery, including the relocation cost from $h^-_t$ to $h^+_{t'}$.
That is, $d_t$ is the direct trip cost, and $d_t + c_{tt'}$ is the cost to serve task $t \in T$ autonomously and relocate to the start of task $t'$.
The autonomous middle miles are discounted by a factor $\alpha \in [0,1]$ and it is assumed that a fraction $\beta \in [0,1)$ of the implied first/last miles are empty.
Additional details on how $\tau_a$ and $c_a$ are calculated are provided by \cite{LeeEtAl2023-OptimizingAutonomousTransfer}.

\begin{figure}[t]
	\begin{mini!}
		%
		{}
		%
		{\sum_{t \in T} d_t + \sum_{a \in A} c_a y_a, \label{eq:athn:obj}}
		%
		{\label{formulation:athn}}
		%
		{}
		%
		%
		\addConstraint
		{\sum_{a \in \delta^+_t} y_a}
		{\le 1 \quad \label{eq:athn:visit}}
		{\forall t \in T,}
		\addConstraint
		{\sum_{a \in \delta^+_t} y_a}
		{= \sum_{a \in \delta^-_t} y_a \label{eq:athn:flowbalance}}
		{\forall t \in T,}
		\addConstraint
		{\sum_{a \in \delta^+_0} y_a}
		{\le K, \label{eq:athn:vehicles}}
		{}
		\addConstraint
		{x_{t'}}
		{\ge x_t + \tau_{tt'} - M (1-y_{tt'})\quad \label{eq:athn:timemtz}}
		{\forall t, t' \in T, (t,t') \in A,}
		\addConstraint
		{x_t}
		{\in [p_t-\Delta, p_t+\Delta] \label{eq:athn:x}}
		{\forall t \in T,}
		\addConstraint
		{y_a}
		{\in \mathbb{B} \label{eq:athn:y}}
		{\forall a \in A.}
	\end{mini!}%
	\vspace{-0.6cm}
	\caption{Mixed Integer Program for Optimizing ATHN Operations.}
\end{figure}

Let $y_a$ be a binary variable that indicates that arc $a\in A$ is selected, and let $x_t$ be the start time of task $t \in T$.
For convenience, let $\delta^+_v$ and $\delta^-_v$ denote the out-arcs and in-arcs of vertex $v \in V$, respectively.
For a given number of autonomous trucks $K$, MIP~\eqref{formulation:athn} asks for a set of at most $K$ routes from source to sink that cover different tasks.
Objective~\eqref{eq:athn:obj} minimizes the total cost.
If a task is not covered it means it is served by conventional means and the cost is that of a direct trip $d_t$.
If the task is covered, the cost differential $c_{tt'}$ is added to calculate the autonomous cost.
Constraints~\eqref{eq:athn:visit} enforce that every task is performed at most once.
Constraints~\eqref{eq:athn:flowbalance} are flow conservation constraints, and Constraint~\eqref{eq:athn:vehicles} limits the number of vehicles to $K$.
Constraints~\eqref{eq:athn:timemtz} are Miller-Tucker-Zemlin constraints \cite{MillerEtAl1960-IntegerProgrammingFormulation} that ensure sufficient time passes between subsequent tasks, where $M$ is a sufficiently large constant.
These constraints also eliminate cycles.
Finally, the variables and their domains are given by Equations~\eqref{eq:athn:x} and \eqref{eq:athn:y}.
The MIP Model~\eqref{formulation:athn} is solved with a blackbox solver after applying the acceleration techniques detailed in \cite{LeeEtAl2023-OptimizingAutonomousTransfer}.
The arc-flows ($y$-variables) are translated into a set of routes, which are represented as sequences of tasks, by tracing the flows from the source node $0$ to the sink node $\lvert T \rvert+1$.
Note that for given optimal routes, the start times ($x$-variables) are typically not unique.
For consistent analysis, the start times are shifted to as early as possible in post processing.

\subsection{Minimizing Required Hub Capacities}
The routes and schedule obtained from the MIP do not take hub utilization into account.
Therefore, it may happen that many trucks are loading and unloading at the same hub at the same time.
To address this issue, a CP model is introduced to shift the schedule to minimize the necessary loading/unloading capacity at the hubs, while satisfying the original time windows and maintaining the same route (sequence of tasks) for every vehicle.

The CP model is based on $\mathcal{K}$; the set of routes obtained from MIP~\eqref{formulation:athn}.
Every route $k \in \mathcal{K}$ is split into an sequence of jobs $\mathcal{J}_k = \{1, 2, \hdots\}$.
For notational convenience, the jobs are numbered sequentially by the order in which they are performed, rather than using the original task numbers.
Figure~\ref{fig:route_jobs} provides a visualization, which will serve as a running example.
Each job $j \in \mathcal{J}_k$ has up to four properties:\\
\begin{itemize}
	\item $\type(j) \in \{\load, \drive, \unload, \relocate, \park\}$; type of job: load, drive, unload, relocate, park, respectively.
	\item $\duration(j)$ (only for types \load, \drive, \unload, \relocate); duration of this job.
	\item $\task(j) \in T$ (only for types \load, \drive, \unload); task associated with this job.
	\item $\hub(j) \in H$ (only for types \load, \unload, \park); hub associated with this job.
\end{itemize}

Each route is modeled with the same repeating sequence of loading at the origin hub (\load), waiting before driving (\park), driving (\drive), waiting at the destination hub before unloading (\park), unloading (\unload), waiting before relocating to the next task if any (\park), relocating (\relocate), and waiting at the origin hub of the next task until loading (\park).
Note that jobs of type \load, \drive, \unload, \relocate{} have a fixed duration, while the duration of \park{} jobs is flexible and can be zero.
When no relocation is necessary (the previous destination is equal to the next origin), the \relocate{} job is still defined with duration zero for convenience.
Every \load, \drive, \unload{} job is trivially associated with an original task $t \in T$.
Jobs \load, \unload, \park{} for which the vehicle is standing still are associated with a hub as described above.

\begin{figure}[t]
	\centering
	\begin{tikzpicture}[node distance=0 cm,outer sep = 0pt]
		
		\tikzstyle{job}=[rectangle,draw, minimum height=1cm, anchor=north west,text centered,text width=2 em]
		
		\tikzstyle{load}=[job, fill=loadcolor]
		\tikzstyle{unload}=[job, fill=unloadcolor]
		\tikzstyle{relocate}=[job, fill=relocatecolor, text=white]
		\tikzstyle{park}=[job, fill=parkcolor]
		\tikzstyle{drive}=[job, fill=drivecolor, text=white]
		
		\node[load] (l1) at (0,0) {\load};
		\node (skip1) [right = of l1] {...};
		\node[relocate] (r1) [right = of skip1] {\relocate};
		\node[park] (r1p) [right = of r1] {\park};
		\node[load] (l2) [right = of r1p] {\load};
		\node[park] (l2p) [right = of l2] {\park};
		\node[drive] (d2) [right = of l2p] {\drive};
		\node[park] (d2p) [right = of d2] {\park};
		\node[unload] (u2) [right = of d2p] {\unload};
		\node[park] (u2p) [right = of u2] {\park};
		\node[relocate] (r2) [right = of u2p] {\relocate};
		\node[park] (r2p) [right = of r2] {\park};
		\node[load] (l3) [right = of r2p] {\load};
		\node (skip3) [right = of l3] {...};
		\node[unload] (u3) [right = of skip3] {\unload};
		
		\draw (l1.south west) -- ($(l1.south west)+(0,-0.5)$) node[right] {$S^k_1$};
		\draw (l2.south west) -- ($(l2.south west)+(0,-0.5)$) node[right] {$S^k_{j-2}$};
		\draw (d2.south west) -- ($(d2.south west)+(0,-0.5)$) node[right] {$S^k_j$};
		
		\draw (u2.south west) -- ($(u2.south west)+(0,-0.5)$) node[right] {$S^k_{j'-2}$};
		\draw (r2.south west) -- ($(r2.south west)+(0,-0.5)$) node[right] {$S^k_{j'}$};
		\draw (l3.south west) -- ($(l3.south west)+(0,-0.5)$) node[right] {$S^k_{j'+2}$};
		
		\draw (u3.south west) -- ($(u3.south west)+(0,-0.5)$) node[right] {$S^k_{\lvert \mathcal{J}_k \rvert}$};
		\draw (u3.south east) -- ($(u3.south east)+(0,-0.5)$) node[right] {$S^k_{\lvert \mathcal{J}_k \rvert + 1}$};
		
		\tikzstyle{line} = [draw, latex'-latex']
		\path [line] ($(r1p.north west)+(0,0.15)$) -- node [pos=0.05,above,align=left] {Hub $h^+_t$} ($(l2p.north east)+(0,0.15)$);
		\path [line] ($(d2p.north west)+(0,0.15)$) -- node [pos=0.95,above,align=right] {Hub $h^-_t$} ($(u2p.north east)+(0,0.15)$);
		\path [line] ($(l2.north west)+(0,0.30)$) -- node [midway,above,align=center] {Task $t \in T$} ($(u2.north east)+(0,0.30)$);
		
	\end{tikzpicture}
	\caption{Jobs $\mathcal{J}_k$ for Route $k \in \mathcal{K}$ (jobs $j, j' \in \mathcal{J}_k$ are used as examples in the main text).}
	\label{fig:route_jobs}
	\begin{mini!}
		%
		{}
		%
		{\sum_{h \in H} C_h, \label{eq:cp:obj}}
		%
		{\label{formulation:cp}}
		%
		{}
		%
		%
		\addConstraint
		{}
		{I^k_j = \texttt{Interval}([S^k_j, S^k_{j+1}])\qquad\qquad \label{eq:cp:intervals}}
		{\forall k \in \mathcal{K}, j \in \mathcal{J}_k,}
		\addConstraint
		{}
		{\mathrlap{\texttt{Cumulative}\left([I^k_j \vert k \in \mathcal{K}, j \in \mathcal{J}_k, \type(j) \in \{\load, \unload\}, \hub(j)=h], C_h\right) \quad \forall h \in H,} \label{eq:cp:cumulative}}
		{}
		\addConstraint
		{}
		{S^k_{j+1} = S^k_j + \duration(j) \label{eq:cp:job_duration}}
		{\forall k \in \mathcal{K}, j \in \mathcal{J}_k, \type(j) \in \{\load, \drive, \unload, \relocate\},}
		\addConstraint
		{}
		{\textrm{dom}(S^k_j) = [\underbar{$s$}^k_j, \overline{s}^k_j] \label{eq:cp:domain_S}}
		{\forall k \in \mathcal{K}, j \in \mathcal{J}_k.}
	\end{mini!}%
	\caption{Constraint Programming Model for Minimizing Required Hub Capacities.}
\end{figure}

\paragraph*{CP Model}
The CP Model~\eqref{formulation:cp} is based on variables $S^k_j$ that indicate the start time of job $j\in \mathcal{J}_k$ in route $k \in \mathcal{K}$, and variables $C_h$ that indicate the necessary loading/unloading capacity at hub $h \in H$.
For convenience, $S^k_{\lvert \mathcal{J}_k \rvert + 1}$ is defined to represent the time at which the final job \unload{} is completed (see Figure~\eqref{fig:route_jobs}).
Objective~\eqref{eq:cp:obj} minimizes the total necessary hub capacity.
To calculate the hub capacity, the model first defines an interval variable $I^k_j$ for every job (Equation~\eqref{eq:cp:intervals}), which implicitly enforces $S^k_j \le S^k_{j+1}$.
Equation~\eqref{eq:cp:cumulative} defines a cumulative constraint for every hub $h\in H$ to collect the intervals of the jobs with type \load{} and \unload{} at that hub, and to assign the neccesary hub capacity to $C_h$.
Note that this cumulative constraint has a variable as the capacity.
It is also worth mentioning that Constraints~\eqref{eq:cp:cumulative} span all the vehicles.
Constraints~\eqref{eq:cp:job_duration} ensure the correct job duration when a duration is defined (\park{} jobs are flexible).
Finally, Equation~\eqref{eq:cp:domain_S} defines the domains of the $S$-variables, where constants $\underbar{$s$}^k_j$, $\overline{s}^k_j$ remain to be defined.
This paper focuses on loading/unloading capacity, but note that CP Model~\eqref{formulation:cp} is easily adapted to other objectives, such as minimizing parking space.

To ensure that the time flexibility $\Delta$ is respected, the domains of the $S$-variables need to be defined accordingly.
For job $j \in \mathcal{J}_k$ of route $k \in \mathcal{K}$ and $\type(j)=\load$, the domain is defined around the desired pickup time of the task to match the MIP: $\textrm{dom}(S^k_j) = [p_{\task(j)} - \Delta, p_{\task(j)} + \Delta]$.
This domain is translated to the following \drive{} and \unload{} jobs to make sure that the flexibility is not exceeded until the task is complete.
With slight abuse of notation, this gives the translated domains $\textrm{dom}(S^k_j) = \textrm{dom}(S^k_{j-2}) + \duration(j-2)$ for jobs of $\type(j) \in \{\drive, \unload\}$ (see Figure~\ref{fig:route_jobs}).

\paragraph*{Redundant Bounds and Constraints}
Non-trivial bounds for jobs of type \relocate{} and \park{} are not strictly necessary, but it is straightforward to derive the following:
\begin{equation}
\textrm{dom}(S^k_{j}) = \begin{cases}
	[\underbar{$s$}^k_{j-2} + \duration(j-2), \overline{s}^k_{j+2} - \duration(j)] & \textrm{ if } \type(j) = \relocate,\\
	[\underbar{$s$}^k_{j-1} + \duration(j-1), \overline{s}^k_{j+1}] & \textrm{ if } \type(j) = \park.
\end{cases}
\end{equation}
Relocation can only start after unloading is completed, and has to start to be in time for the next loading (see job $j'$ in Figure~\ref{fig:route_jobs}).
In a similar way, parking can only start after the previous job is completed, and parking for duration zero is possible until the upper bound of the next job.

For the current Objective~\eqref{eq:cp:obj}, the following observation is used to reduce the search space: To minimize the loading/unloading capacity, it does not matter when relocation takes place between the jobs \unload{} and \load{}.
It is therefore possible without loss of generality to impose that relocation starts immediately after unloading:
\begin{equation}
	S^k_{j} = S^k_{j-1} \quad \forall k \in \mathcal{K}, j \in \mathcal{J}_k, \type(j) = \relocate.
\end{equation}

%% file: sections/case_study.tex
\section{Case Study}
\label{sec:case_study}
This paper conducts a realistic case study based on data from Ryder System, Inc. (Ryder), one of the largest transportation and logistics companies in North America.
Ryder has provided a dataset that is representative for its dedicated transportation business in the US, reducing the scope to orders that are strong candidates for automation.
Following \cite{LeeEtAl2023-OptimizingAutonomousTransfer}, the case study focuses on orders that are \emph{challenging} in the sense that they would currently induce an empty return trip.
Every order represents a load with an origin, a destination, and a scheduled release time.
The hubs are chosen based on data from October to December 2019, while the experiments are based on 6842 loads in the first four weeks of October.

\paragraph*{Experimental Settings}

\begin{figure*}[t]
	\scriptsize
	\centering
	\begin{tabular}{ll}
		\toprule
		Parameter & Value \\
		\midrule
		$n$ & 6842 loads\\
		$\lvert H \rvert$ & 100 transfer hubs\\
		$\gamma$ & 40\% discount for autonomous mileage during hub-assignment\\
		$\beta$ & 25\% first/last-mile inefficiency\\
		$\alpha$ & 25\% discount for autonomous mileage\\
		$\Delta$ & 1 hour pickup-time flexibility\\
		$\sigma$ & 30 minutes autonomous truck loading/unloading time\\
		$K$ & 100 autonomous trucks\\
		\bottomrule
	\end{tabular}
	\vspace{0.2cm}
 	\captionof{table}{Baseline Parameter Values for the Case Study.}%
	\label{tab:base_parameters}%
\end{figure*}

The parameter settings are taken from \cite{LeeEtAl2023-OptimizingAutonomousTransfer} and summarized in Table~\ref{tab:base_parameters}.
All methods from Section~\ref{sec:methodology} are implemented in Python 3.9.
The MIP models are solved with Gurobi 9.5.2, and the CP models are solved with CP-SAT 9.6~\cite{ortools}.
All time-related data is rounded to the nearest minute to facilitate the integer domains that are required by CP-SAT.
The experiments are conducted on a Linux machine with dual Intel Xeon Gold 6226 CPUs on the PACE Pheonix cluster~\cite{PACE2017-PartnershipAdvancedComputing}.
Each experiment is assigned to use at most 24 cores and 192GB of RAM.
If Gurobi runs out of memory, the experiment is repeated on a machine with 384GB of RAM, and if that fails, the number of cores is halved until the solver terminates successfully.
CP-SAT did not encounter memory issues at any point.
The Gurobi time limit is set to three hours, except for the sensitivity analysis for flexibility $\Delta$, which is given 12 hours to obtain better solutions and is warm started with a MIP start.
The CP-SAT solver is warm started with the solution obtained by the MIP, which is provided as a hint.

%% file: sections/results.tex
\section{Results}
\label{sec:results}

\begin{figure}[!t]
	\centering
 	\begin{minipage}{0.5\textwidth}
		\centering
		\includegraphics[height=0.65\textwidth]{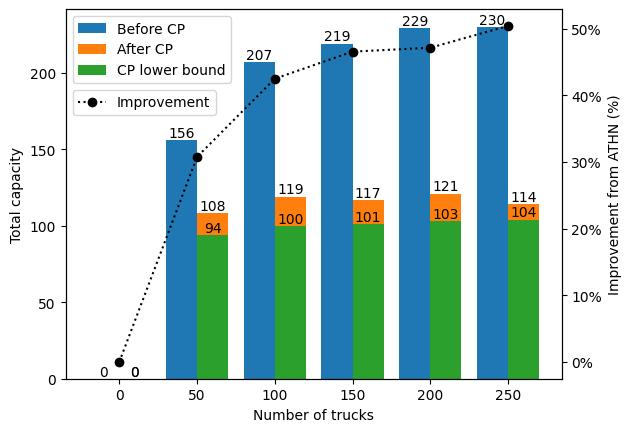}
		\caption{Capacity Before and After CP.}
		\label{fig:base_cp_results}
	\end{minipage}%
	\begin{minipage}{0.5\textwidth}
		\centering
		\includegraphics[height=0.65\textwidth]{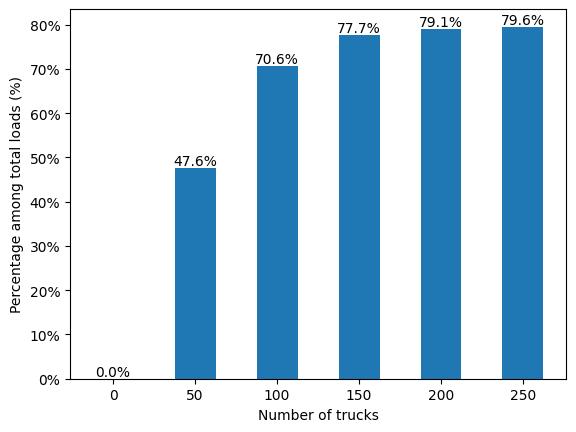}
		\caption{Loads Served Autonomously.}
		\label{fig:base_loads_served}
	\end{minipage}
\end{figure}

\begin{figure}[!t]
	\centering
	\begin{subfigure}{0.47\textwidth}
		\centering
		\includegraphics[width=\textwidth]{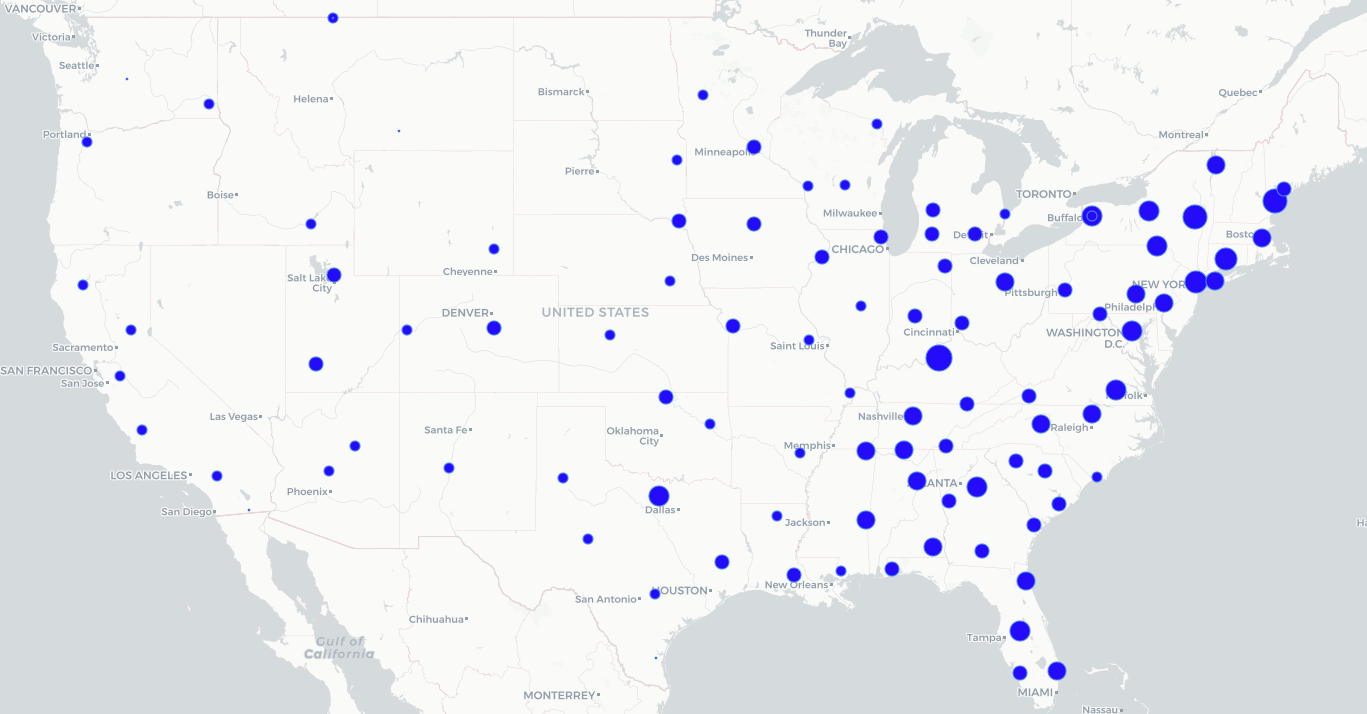}
		\caption{Before CP Optimization}
		\label{fig:before_cp_map}
	\end{subfigure}
	\begin{subfigure}{0.47\textwidth}
		\centering
		\includegraphics[width=\textwidth]{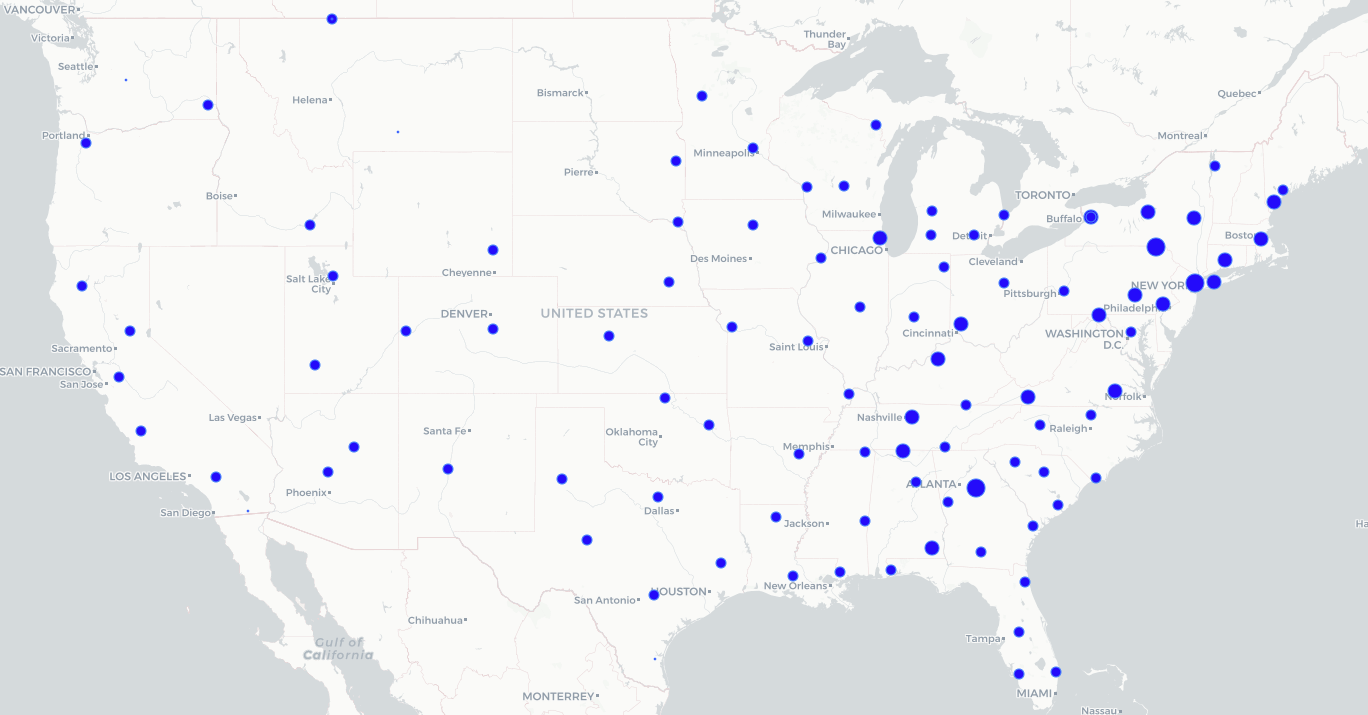}
		\caption{After CP Optimization}
		\label{fig:after_cp_map}
	\end{subfigure}
	\caption{Capacity Reduction for the 100 Truck Case (circle area proportional to hub capacity).}
	\label{fig:before_after_map}
\end{figure}

Figure~\ref{fig:base_cp_results} summarizes the results of using the CP model to minimize the necessary hub capacity for loading and unloading trucks.
``Before CP'' shows the required capacity if the MIP solution were implemented immediately, while ``After CP'' shows the results after applying the CP model.
The ATHN problem was solved to optimality for all instances except when $K=50$, which remained at a 0.05\% optimality gap within the given time limit.
For the CP problem, all instances were solved to optimality in under 30 seconds. 
The figure shows that the CP model is highly effective in reducing the total hub capacity for all instances, reducing the necessary capacity by 31\% up to 50\%.
For the base case of $K=100$ the CP model can reduce capacity by 42\%.
Note that this may correspond to significant monetary savings: If each unit of loading/unloading capacity requires the assistance of a mechanic around the clock (three shifts of \$57,557 per year \cite{Glassdoor}), the CP model reduces the annual labor cost by \$15.2 million.
Figure~\ref{fig:base_cp_results} also shows that the resulting capacity is close to the lower bound \emph{for any ATHN solution}, which is obtained by removing the time constraints between subsequent tasks, i.e., treat every task as if it is the only task on the route.

As the number of trucks increases, the ATHN starts serving more loads autonomously, as shown by Figure~\ref{fig:base_loads_served}.
Without the CP model, this leads to a substantial increase in hub capacity as more loads are added to the system.
The CP model completely mitigates this effect, and is able to maintain an almost stable hub capacity as the workload increases.
This reveals a suprising robustness to accomodate new orders that the CP model is able to exploit.
When making investment decisions and hiring personnel to operate the hubs, this is a very desirable property.
The maps in Figure~\ref{fig:before_after_map} demonstrate that the necessary capacity is reduced throughout the system, and the peaks in the South and Northeast have decreased significantly.
The capacity at the hub near Louisville, Kentucky for example was reduced from 7 to 2.

\begin{figure}[!t]
	\centering
	\begin{subfigure}{0.47\textwidth}
		\centering
		\includegraphics[width=\textwidth]{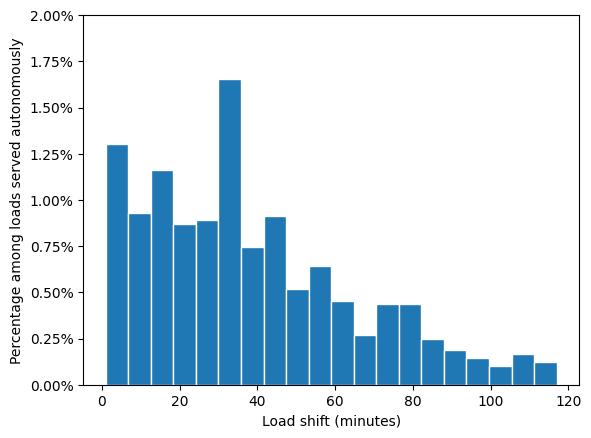}
		\caption{Histogram of Shifted Loading Times.}
		\label{fig:base_positive_load_shift}
	\end{subfigure}
	\begin{subfigure}{0.47\textwidth}
		\centering
		\includegraphics[width=\textwidth]{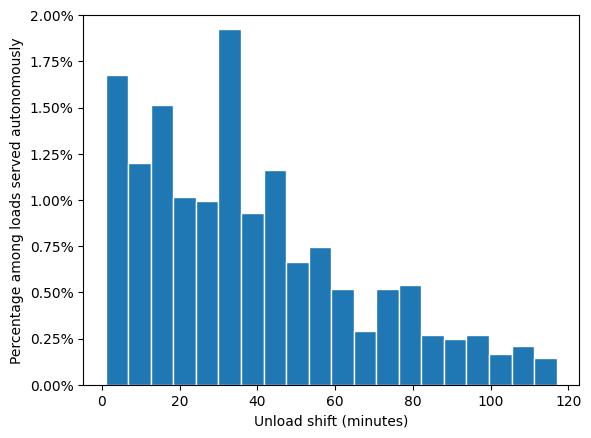}
		\caption{Histogram of Shifted Unloading Times.}
		\label{fig:base_positive_unload_shift}
	\end{subfigure}
	\caption{Shifted Schedule Times for the 100 Truck Case.}
	\label{fig:base_shift}
\end{figure}

Interestingly, only a small portion of the jobs needs to be rescheduled to obtain the substantial savings in hub capacity.
E.g., in the 100-truck case only 13\% of loadings were rescheduled.
Figure~\ref{fig:base_shift} provides histograms for the absolute size of the shift for the loading and unloading times that were shifted.
It can be seen that a majority of these jobs were moved by no more than 30 minutes, which shows that the ATHN solutions are sufficiently flexible to be adjusted without propagating delays through the schedule.
The fact that the CP model is easy to solve and only makes small modifications also makes it an attractive tool at the operational level: if order details are changed or delays are encountered, the CP model can quickly be re-solved to avoid causing overlap in loading and unloading at the hubs.

\begin{figure}[!t]
	\centering
	\begin{subfigure}{0.47\textwidth}
		\centering
		\includegraphics[width=\textwidth]{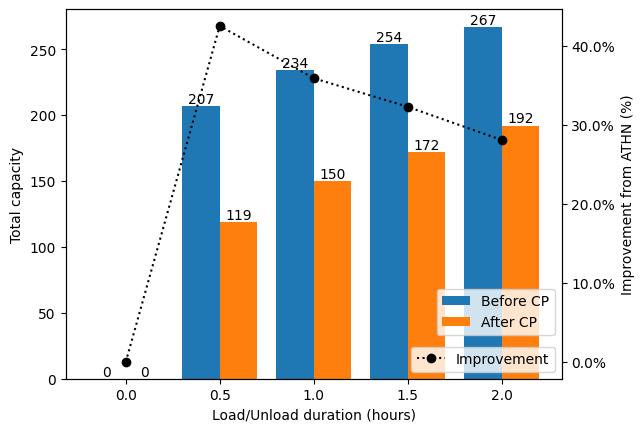}
    	\caption{Impact of Loading/Unloading Time $\sigma$.}
    	\label{fig:load_unload}
	\end{subfigure}
	\begin{subfigure}{0.47\textwidth}
		\centering
		\includegraphics[width=\textwidth]{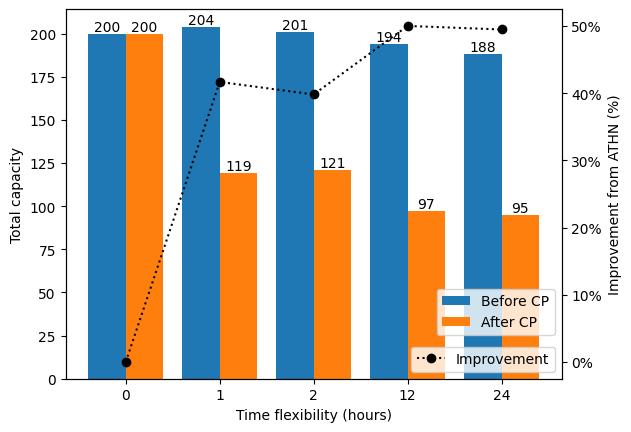}
		\caption{Impact of Flexibility $\Delta$.}
		\label{fig:time_flexibility}
	\end{subfigure}
	\caption{Sensitivity Analysis Loading/Unloading Time and Flexibility for the 100 Truck Case.}
	\label{fig:parameters}
\end{figure}

\paragraph*{Sensitivity Analysis}
The loading/unloading time $\sigma$ and the flexibility $\Delta$ are two significant factors that affect the ATHN.
Figure~\ref{fig:parameters} summarizes how these parameters impact the total hub capacity and the CP model's performance.
The ATHN is solved for loading/unloading time from zero to two hours, and for no flexibility up to one day of flexibility.
Note that loading/unloading time up to two hours may be realistic if the vehicle is inspected every time it leaves or enters a hub.

For $\sigma$, the ATHN solver found optimal solutions for all instances and the CP problem found optimal schedules within 15 seconds.
As $\sigma$ increases, the total hub capacity increases because jobs overlap more frequently.
The CP model is not able to fully compensate for this effect, but can still improve the hub capacity by at least 28\%.

For $\Delta$, optimal solutions were found for $\Delta \in \{0\textrm{h}, 1\textrm{h}\}$, a near-optimal solution within 0.05\% optimality gap was found for $\Delta = 2\textrm{h}$, while $\Delta \in \{12\textrm{h}, 24\textrm{h}\}$ remained at a 3\% optimality gap.
The CP solver took 7 minutes for $\Delta =12\textrm{h}$ and 13 minutes for $\Delta = 24\textrm{h}$ to find optimal schedules, which suggests that the CP computation time is more sensitive to changes in $\Delta$ than in $\sigma$.
This is explained by the fact that increasing the flexibility significantly increases the search space.
Again, the CP model shows substantial improvement over the initial solution with improvements ranging from 39\% to 50\% for the case study.
When the flexibility increases, the necessary capacity goes down even before CP is applied, but the CP model is able to benefit more from the additional freedom.

\paragraph*{Impact of Redundant Bounds and Constraints}
Table~\ref{tab:speedups} shows the impact of the redundant bounds and constraints introduced in Section~\ref{sec:methodology}.
These results are based on the baseline parameter settings for different numbers of trucks.
It can be seen that only adding redundant bounds can both speed up or slow down the solver.
For example, the solver becomes 1.4 times faster for $K=100$ trucks, but a factor 0.8 slower for $K=250$ trucks.
The main benefit comes from adding redundant constraints, which speeds up the solver by up to 7.4 times for $K=150$ trucks.
When the redundant bounds and constraints are combined, performance may be improved further, as is the case for $K=250$ trucks, but it appears that only including redundant constraints is the most efficient setting for these experiments.

\begin{figure*}[!t]
	\renewcommand{\arraystretch}{1.2}
	\scriptsize
	\centering
	\begin{tabular}{ccccccc}
		\toprule
		& & \multicolumn{5}{c}{Number of Trucks} \\
		\cmidrule(lr){3-7}
		Redundant Bounds & Redundant Constraints & 50    & 100   & 150   & 200   & 250 \\
		\midrule
		\xmark & \xmark & 1.0x  & 1.0x  & 1.0x  & 1.0x  & 1.0x \\
		\cmark & \xmark & 1.1x  & 1.4x  & 1.2x  & 0.9x  & 0.8x \\
		\xmark & \cmark & 3.9x  & 6.8x  & 7.4x  & 3.9x  & 5.4x \\
		\cmark & \cmark & 3.5x  & 5.5x  & 6.3x  & 3.1x  & 5.5x \\
		\bottomrule
	\end{tabular}%
	\vspace{0.2cm}
	\captionof{table}{CP Solving Time Speedup compared to Omitting Redundant Bounds and Constraints.}%
	\label{tab:speedups}%
\end{figure*}

%% file: sections/conclusion.tex
\section{Conclusion}
\label{sec:conclusion}

The Autonomous Transfer Hub Network (ATHN) is one of the most promising ways to adapt self-driving trucks for the freight industry.
This paper proposes a framework for optimizing ATHN operations with respect to hub utilization.
To accomplish this, a MIP model generates autonomous vehicle routes and an initial schedule, which is then optimized using a CP model to reduce the required hub capacity for loading and unloading trucks. Results from the Ryder case study demonstrate the effectiveness of this approach, with the CP model reducing the total capacity by at least 42\% while requiring only minor schedule modifications.
This may save \$15.2M per year in labor cost and is close to the lowest possible capacity for any initial schedule.
Sensitivity analysis on loading/unloading duration and flexibility provides practical insights for ATHN operations and systematically shows the benefit of the CP model.
Future work may explore alternative objective functions, such as minimizing parking space.
Another interesting direction is to attempt to jointly optimize routes and hub capacity directly.